\documentclass[submission,copyright,creativecommons]{eptcs}
\pdfoutput=1
\providecommand{\event}{Preprint} %
\usepackage{iftex}
\usepackage{amsmath}
\usepackage{amssymb}
\usepackage{placeins}
\usepackage{physics}
\usepackage{placeins}
\usepackage{floatrow}                             %
\usepackage{framed}
\usepackage{caption}
\usepackage{geometry}
\usepackage{lipsum}
\usepackage{graphicx}
\usepackage{floatrow}
\usepackage{hyperref}

\ifpdf
  \usepackage{underscore}         %
  \usepackage[T1]{fontenc}        %
\else
  \usepackage{breakurl}           %
\fi

\def\scalemed{0.55}
\def\scalebig{0.70}

\newtheorem{theorem}{Theorem}
\title{Functor String Diagrams: A Novel Approach to Flexible Diagrams for Applied Category Theory}
\author{Vincent Abbott
\institute{College of Engineering, Computing and Cybernetics\\
Australian National University\\
Canberra, ACT, Australia}
\email{Vincent.Abbott@anu.edu.au}
\and
Gioele Zardini
\institute{
Laboratory for Information and Decision Systems\\
Massachusetts Institute of Technology\\
Cambridge, MA, USA}
\email{gzardini@mit.edu}
}
\def\titlerunning{Functor String Diagrams}
\def\authorrunning{V. Abbott, G. Zardini}

\begin{document}
\maketitle

\begin{abstract}
The study of abstraction and composition---the focus of category theory---naturally leads to sophisticated diagrams which can encode complex algebraic semantics. 
Consequently, these diagrams facilitate a clearer visual comprehension of diverse theoretical and applied systems.
Complex algebraic structures---otherwise represented by a forest of symbols---can be encoded into diagrams with intuitive graphical rules. 
The prevailing paradigm for diagrammatic category theory are monoidal string diagrams whose specification reflects the axioms of monoidal categories. 
However, such diagrams struggle in accurately portraying crucial categorical constructs such as functors or natural transformations, obscuring central concepts such as the Yoneda lemma or the simultaneous consideration of hom-functors and products.
In this work, we introduce \emph{functor string diagrams}, a systematic approach for the development of categorical diagrams which allows functors, natural transformations, and products to be clearly represented. 
We validate their practicality in multiple dimensions. 
We show their usefulness for theoretical manipulations by proving the Yoneda lemma, show that they encompass monoidal string diagrams and hence their helpful properties, and end by showing their exceptional applied utility by leveraging them to underpin neural circuit diagrams, a method which, at last, allows deep learning architectures to be comprehensively and rigorously expressed.
\end{abstract}

\section{Introduction}
Category theory enables a rigorous and abstract understanding of composition, offering tools for studying a wide range of systems.
One of such tools are diagrams. 
Cohesive diagrams need to obey strict compositionality rules, making category theory the natural means by which to approach them. 
In describing systems, diagrams offer a rigorous framework to visualize and analyze them.
Furthermore, by diagramming a system, its fundamental abstract structure is revealed.

Categorical diagrams encode the compositional structure of systems, naturally reflecting their algebraic semantics and exchanging the management of symbols for the manipulation of diagrams.
The most prominent diagrams in category theory are monoidal string diagrams \cite{selinger-survey-2009, piedeleu-introduction-2023}.

These have been applied, or found to underlie, diagrams in various domains. This includes Feynman diagrams in physics, Petri nets for engineered systems, diagrams for quantum computing, and shared theoretical structures in physics, logic, and topology \cite{abramsky-categorical-2008, baez-physics-2010, meseguer-petri-1990}.
Monoidal string diagrams offer a common unified approach, letting theorems be applied across domains.

Current diagrammatic approaches, however, are specialized to focus on certain algebraic concepts. They employ topological isotopy rules which ensure that diagrams correspond to topological manipulations which shift around components. 
This limits them. Monoidal string diagrams, for instance, rely on a privileged product to combine morphisms. This struggles to represent systems where the actions of both linear and non-linear operations are critical to represent, for instance. Furthermore, monoidal string diagrams lack a natural expression for functors and, hence, fail to express many insightful categorical concepts.

Alternative string diagrams \cite{marsden-category-2014, nakahira-diagrammatic-2023} address functors and natural transformations but lack the means to readily express products. These diagrams use colors to represent separate categories, and strings are functors showing boundaries between them. Natural transformations are nodes on these boundaries, and may introduce additional functors, acting as intersections between boundaries. These diagrams are focused on advanced categorical concepts, the behaviour of monads and adjunctions, for instance. A variety of existing diagrammatic approaches are showcased in the appendix, Section \ref{sec:existing}.

For applied cases, however, streamlined diagrams are preferred. While morphisms, functors, and natural transformations frequently arise in practice, color coding categories impose distraction but little benefit. Furthermore, expressing both functors and products in either scheme is unclear and obfuscates the underlying expression. Proofs related to areas in which diagrams are not-specialized are difficult to parse. The topological formality of current diagrams, then, can be seen as a hindrance to generalization, cross-domain understanding, and application.

Instead, it will be valuable to distill the construction of diagrams down to clear principles which ensures we can construct a reliable correspondence between categorical expressions and diagrams. Rather than specializing diagrams for a particular case, we will ensure they can robustly represent categorical expressions and will introduce new features as needed. These new features will each provide local graphical intuition and will correspond with algebraic rules, but will not be required to be consistent with an overarching topological regularity. We will prefer to use colors for annotation rather than expression, granting additional utility for explaining and applying diagrams. 

This leads to functor string diagrams \cite{abbott-robust-2023}, a principled approach to constructing diagrams for various systems. We will demand that diagrams can be readily corresponded to standard object-morphism-object-etc. expressions. New features will be introduced as needed, with each implementing an algebraic rule in a graphically intuitive manner that motivates its use. Therefore, they can reveal the algebraic structure in a range of cases, from formal proofs to the design of deep learning models, in a clearly understandable manner.

\subsection*{Our Contributions} In this work, the two principles underlying functor string diagrams, vertical section decomposition and equivalent expression for graphical intuition, are introduced. 
The disciplined diagrams constructed from these principles are shown to provide both theoretical and applied utility. 
We show this utility by presenting a straightforward proof of the Yoneda lemma, showing existing monoidal string diagrams to be a specific instance of functor string diagrams, and finally show that they provide mathematical rigor to neural circuit diagrams, addressing the open-problem of comprehensively communicating the details of deep learning models.

\section{The Foundations of Functor String Diagrams}\label{sec:foundations}
Functor string diagrams fundamentally take the perspective that categorical diagrams should be flexible yet rigorous, adaptable to various use cases while maintaining clear correspondence to underlying symbolic expressions and having features which guide intuition regarding the application of theorems. 
This is achieved by establishing two principles for the construction of diagrams.

This section will introduce and highlight these principles. If they are followed, we will attain flexible, disciplined diagrams that can represent just about any categorical concept. Compared to monoidal string diagrams which strictly demand isotopy, we have the flexibility to introduce functors and natural transformations in a rigorous, intuitive manner. Furthermore, these diagrams can be adapted to different use cases. As long as we are aware of the equivalent expressions employed, we always have diagrams corresponding to unique up to isomorphism categorical expressions.

\subsection{Vertical Section Decomposition}
Functor string diagrams can be decomposed into vertical columns. 
Vertical columns either represent the current object, or a morphism acting on that object. 
Below in Figure \ref{fig:fsd}, we can see a basic functor string diagram mapping through $\displaystyle f:x\rightarrow y$ and $\displaystyle g:x\rightarrow z$ to construct $\displaystyle f;g:x\rightarrow z$. \textit{(}$\displaystyle ;$\textit{ representing forward composition.) }Note how we have vertical columns for the objects $\displaystyle x,y,z$, and other vertical columns for the morphisms $\displaystyle f$ and $\displaystyle g$.
\FloatBarrier

\begin{figure}[h]
  \floatbox[{\capbeside\thisfloatsetup{capbesideposition={left,top}}}]{figure}[\FBwidth]
  {\caption{
		Here, we have a morphism with multiple functor wires present. Using the principle of equivalent expression, this can be re-expressed into an easier-to-understand form.
	}
	\label{fig:fsd}}
  {\includegraphics[scale=\scalemed]{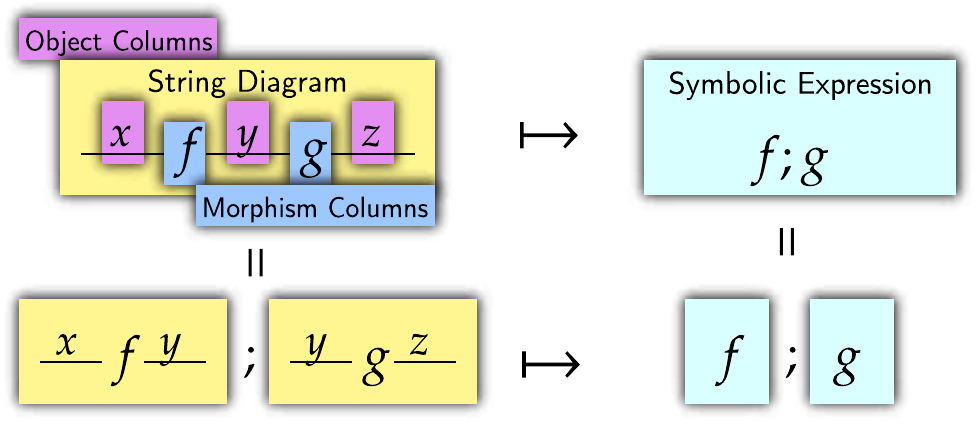}}
\end{figure}

This gives us the first principle for functor string diagrams, vertical section composition, which states that all vertical sections of a diagram correspond to either an object or a morphism. Furthermore, as we express identities by straight lines, the same as objects, the construction of diagrams via vertical section decomposition corresponds with the composition and identity axioms of a category. Therefore, there is a one-to-one correspondence between valid composed functor string diagrams and valid composed symbolic categorical expressions.

\subsection{Equivalent Expression}
In practice, we will want to introduce additional features to diagrams to represent functors, natural transformations, and other constructs. 
To maintain cohesion, we will demand that our means of representing these additional features ultimately corresponds to columns of objects and morphisms, abiding by the first principle. 
This grants the second principle, that all additional features are equivalent expressions of existing features to support graphical intuition.

\begin{figure}[!htb]
    \centering
    \includegraphics[scale=\scalemed]{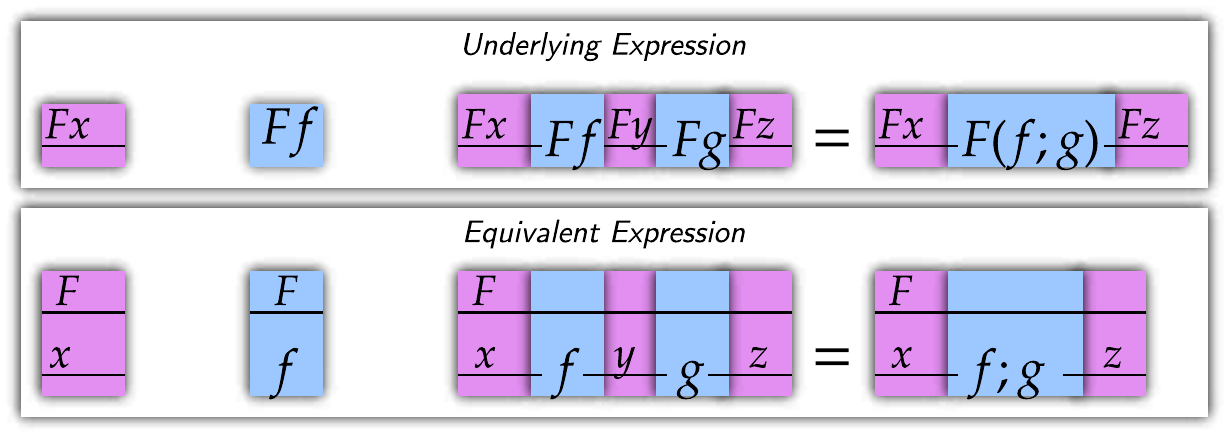}
    \caption{Objects and morphisms in the image of functors can be explicitly expressed. Functors obey the axioms that $f:x\rightarrow y$ maps to $Ff:Fx\rightarrow Fy$ and $F(f;g)=Ff;Fg$. Therefore, we develop an equivalent expression which graphically enforces these axioms, necessarily mapping a morphism along with its objects and necessarily mapping composed morphisms.}
    \label{fig:EE}
\end{figure}

This guides us towards expressing functors. 
Functors are morphisms from~$\displaystyle \mathbf{C}$ to $\displaystyle \mathbf{D}$ which map objects and morphisms so that composition is preserved. 
Fundamentally, we require that~$\displaystyle f:x\rightarrow y$ maps to $\displaystyle Ff:Fx\rightarrow Fy$ and that $\displaystyle F( f;g) =Ff;Fg$. We express an object $\displaystyle Fx$ by a vertical column with a functor wire $\displaystyle F$ placed above $\displaystyle x$, and $\displaystyle Ff$ by a vertical column with the functor wire $\displaystyle F$ placed above $\displaystyle f$, as shown in Figure \ref{fig:EE}. When constructing $\displaystyle Ff;Fg$, we see that graphical intuition supports its equality to $\displaystyle F( f;g)$.

Similarly, we find that natural transformations also have an intuitive equivalent expression. We can write the defining axiom of natural transformations at the top, wherein the functors can be equivalently expressed with functor wires. We notice that the natural transformation can move along the functor line, hence, we establish an equivalent expression for natural transformation components by writing them on the functor wire. This is similar to the methods used by Marsden and Nakahira \cite{marsden-category-2014, nakahira-diagrammatic-2023} and ensures that the axioms of diagrams reflect the axioms of symbols.

\begin{figure}[h]  \floatbox[{\capbeside\thisfloatsetup{capbesideposition={left,top}}}]{figure}[\FBwidth]
  {\caption{
		We develop equivalent expressions to consider natural transformations. We see that natural transformations move along the functor wire and, hence, we show them as morphisms on the functor wire.
	}
	}
  {\includegraphics[scale=\scalemed]{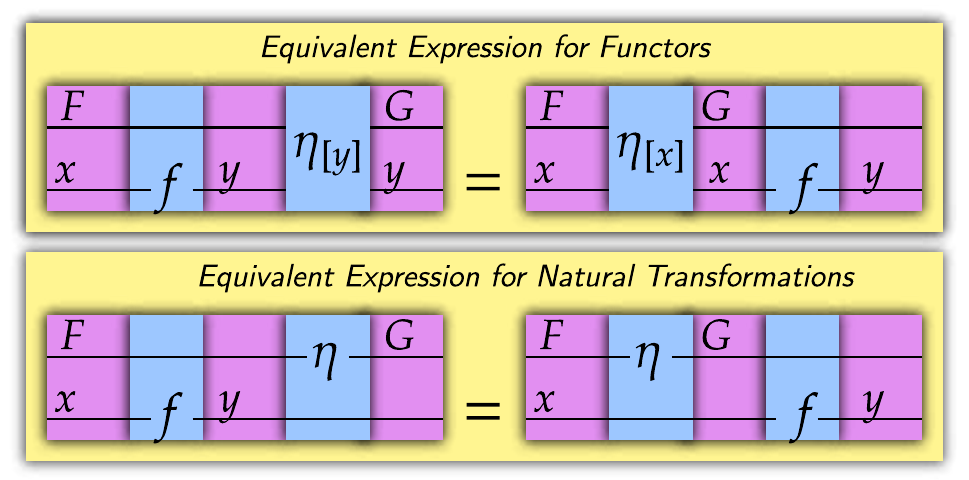}}
\end{figure}

\section{The Yoneda Lemma with Functor String Diagrams} \label{sec:Yoneda}

To prove the algebraic utility of functor string diagrams, we will be expressing the Yoneda lemma in a natural manner.
The Yoneda lemma, a fundamental concept in category theory, has its relevance to various domains potentially obscured when the role of natural transformations is not evident.
As traditional monoidal string diagrams can not consider functors or natural transformations, a novel diagrammatic method is needed to understand it. Here, we will provide a proof of the Yoneda lemma, borrowing from Nakahira, but which uses the principles laid out to offer a brief, intuitive proof.

The Yoneda lemma relates to hom-functors and the category $\displaystyle \mathbf{Set}$. In $\displaystyle \mathbf{Set}$, objects are sets and morphisms are functions between them. In $\displaystyle \mathbf{Set}$, the elements of sets correspond to morphisms from the set of one element, $\displaystyle 1$, to the object corresponding to the set. In $\displaystyle \mathbf{Set}$, morphisms are functions. Functions $f:x\rightarrow y$ are defined by the output element $\varphi ;f:1\rightarrow y$ produced for each input element $\varphi: 1 \rightarrow x$, and any function pairing each input to any output can be constructed.

\begin{figure}[!htb]
  \floatbox[{\capbeside\thisfloatsetup{capbesideposition={left,top}}}]{figure}[\FBwidth]
  {\caption{
		The hom-functor maps morphisms $g:y \rightarrow z$ to a function between sets $\mathbf{C}(x,g):\mathbf{C}(x,y)\rightarrow \mathbf{C}(x,z)$. This function implements composition with $h:x \rightarrow y$ when provided with the element of the set $\mathbf{C}(x,y)$ corresponding to $h$. We can introduce a number of equivalent expressions to more clearly visualize these expressions. This includes not drawing $\displaystyle 1$ object wires.
	}
	\label{fig:hom-ee}}
  {\includegraphics[scale=\scalemed]{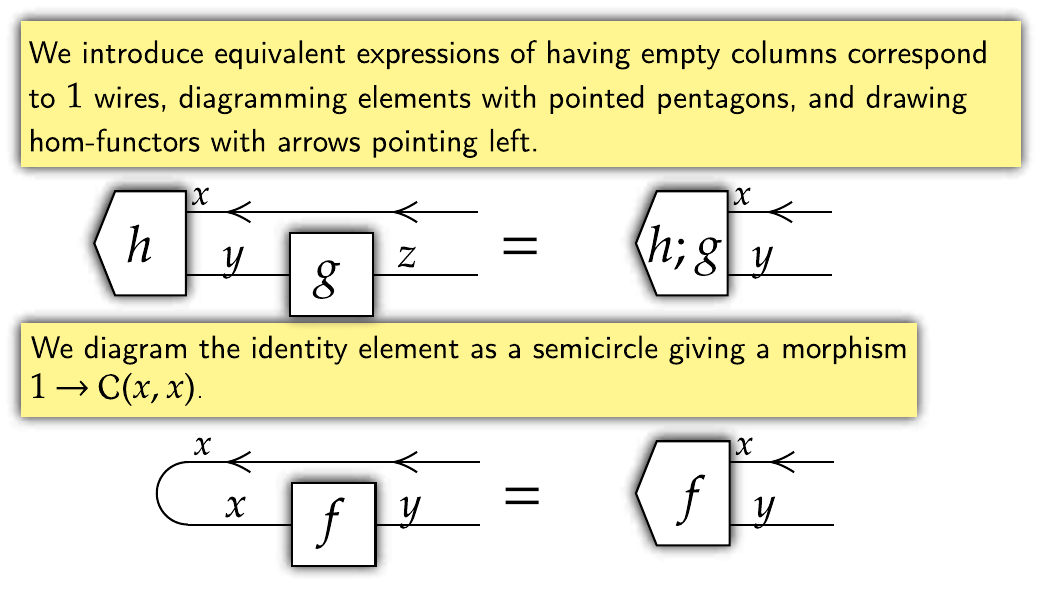}}
\end{figure}

For any category $\displaystyle \mathbf{C}$ where morphisms between objects forms a set \textit{(locally small)}, there is a hom-functor for each object $\displaystyle x$ which maps from $\displaystyle \mathbf{C}$ to $\displaystyle \mathbf{Set}$, sending objects $\displaystyle y$ to the \textit{set} of morphisms $\displaystyle \mathbf{C}( x,y)$ and morphisms $\displaystyle g:y\rightarrow z$ to a morphism $\displaystyle \mathbf{C}( x,g) :\mathbf{C}( x,y)\rightarrow \mathbf{C}( x,z)$. The hom-functored morphism corresponds to the following operation in Figure \ref{fig:hom-ee}.

\FloatBarrier
Regarding hom-functors, the Yoneda lemma makes the following claim.
\begin{theorem}[Yoneda Lemma] The natural transformations $\displaystyle \eta :\mathbf{C}( x,\_)\rightarrow \Phi $, both being functors $\displaystyle \mathbf{C}\rightarrow \mathbf{Set}$, are in one-to-one correspondence with the elements $\displaystyle \varphi \in \Phi x$, the set the object $\displaystyle x\in \text{Ob}(\mathbf{C})$ is mapped to by $\displaystyle \Phi $:
\begin{equation*}
\text{Nat}(\mathbf{C}( x,\_) ,\Phi ) \cong \Phi x.
\end{equation*}
\end{theorem}

\textit{Proof.} We start with an element $\displaystyle \varphi \in \Phi x$. We represent this as a morphism $\displaystyle 1\rightarrow \Phi x$. Diagrammatically, we see that this provides a hook for a morphism $\displaystyle f:x\rightarrow y$ to attach, giving an element $\displaystyle ( \Phi f)( \varphi ) \in \Phi y$. We will define a function mapping from any morphism $\displaystyle f:x\rightarrow y$ to the outcome of this operation, hence, its a morphism $\displaystyle \varphi _{[ y]}^{*} :\mathbf{C}( x,y)\rightarrow \Phi y$. This morphism $\displaystyle \mathbf{C}( x,y)\rightarrow \Phi y$ exists by virtue of the fact that a function mapping from any set of inputs to outputs is allowed in $\displaystyle \mathbf{Set}$.
\FloatBarrier
\begin{figure}[!htb]
    \centering
    \includegraphics[scale=\scalemed]{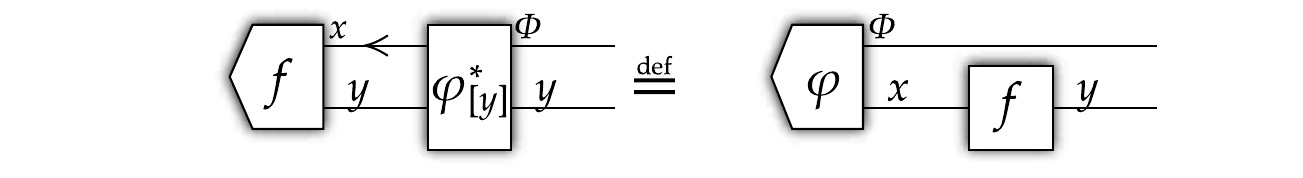}
\end{figure}
We need to confirm whether this is a natural transformation. We do this by testing the above with a generic function $\displaystyle g:y\rightarrow z$.

\FloatBarrier
\begin{figure}[!htb]
    \centering
    \includegraphics[scale=\scalemed]{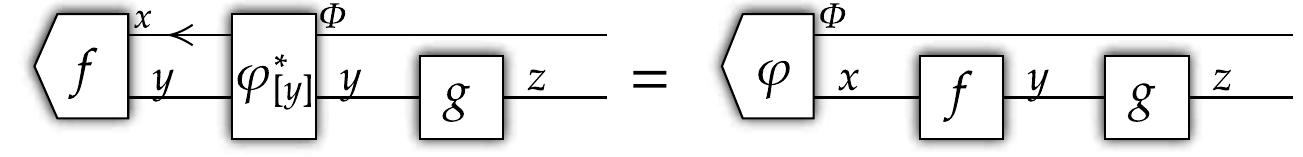}
    \includegraphics[scale=\scalemed]{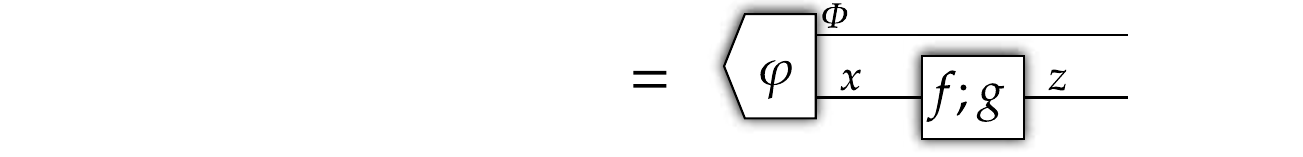}
    \includegraphics[scale=\scalemed]{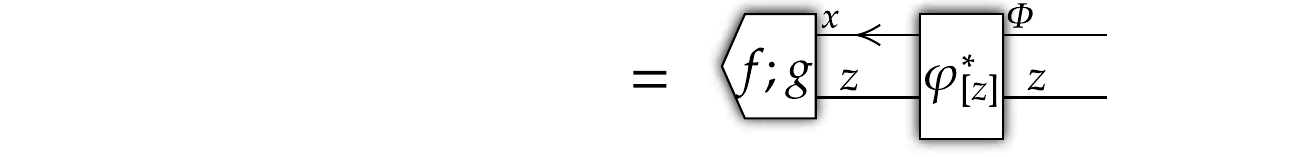}
    \includegraphics[scale=\scalemed]{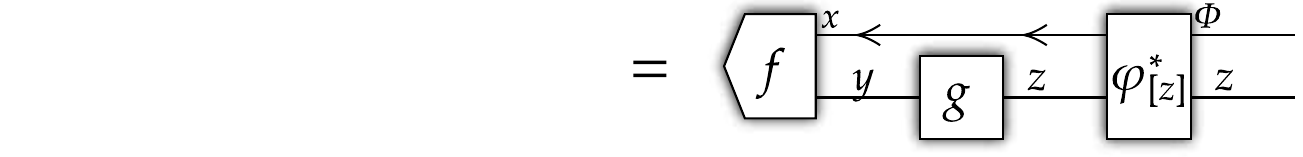}
\end{figure}
\FloatBarrier
These functions are equal for all elements in the domain $\mathbf{C}(x,y)$, making them equal.
\FloatBarrier
\begin{figure}[!htb]
    \centering
    \includegraphics[scale=\scalemed]{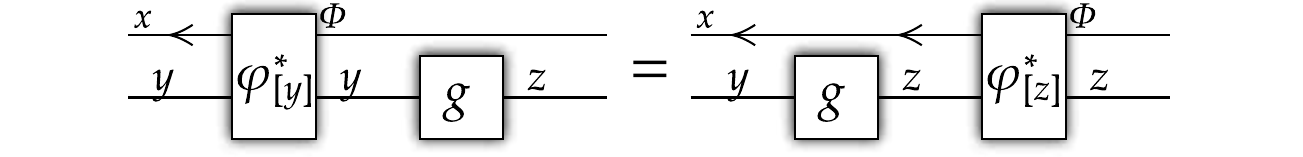}
\end{figure}
\FloatBarrier
\begin{figure}[!htb]
    \centering
    \includegraphics[scale=\scalemed]{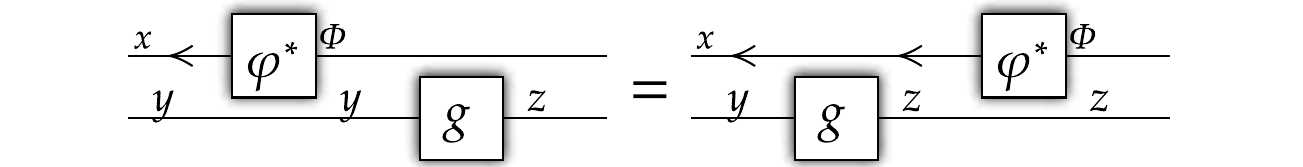}
\end{figure}
\FloatBarrier
Therefore, we have a procedure ``$\displaystyle *$'' which maps from elements $\displaystyle \varphi \in \Phi x$ to natural transformations $\displaystyle \varphi ^{*} :\mathbf{C}( x,\_)\rightarrow \Phi $. The next phase of the proof is to show that there are elements corresponding to natural transformations, and that taken together, these processes form a bijection.

Starting with a natural transformation $\displaystyle \eta :\mathbf{C}( x,\_)\rightarrow \Phi $, we use the identity to obtain an element $\displaystyle 1\rightarrow \Phi x$. This gives a map ``$\displaystyle \dagger $'' from natural transforms to elements, $\displaystyle \dagger :\text{Nat}(\mathbf{C}( x,\_) ,\Phi )\rightarrow \Phi x$.

\FloatBarrier
\begin{figure}[h]
    \centering
    \includegraphics[scale=\scalemed]{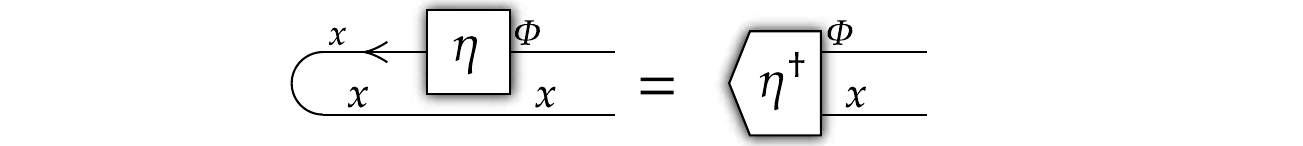}
\end{figure}
\FloatBarrier

Now, we will find the natural transformation $\displaystyle \eta ^{\dagger *}$ corresponding to the element $\displaystyle \eta ^{\dagger }$.
\FloatBarrier
\begin{figure}[!htb]
    \centering
    \includegraphics[scale=\scalemed]{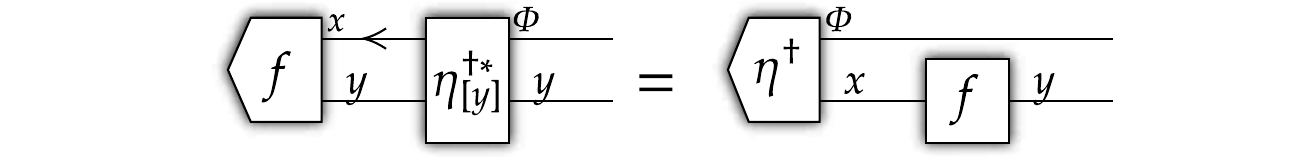}
    \includegraphics[scale=\scalemed]{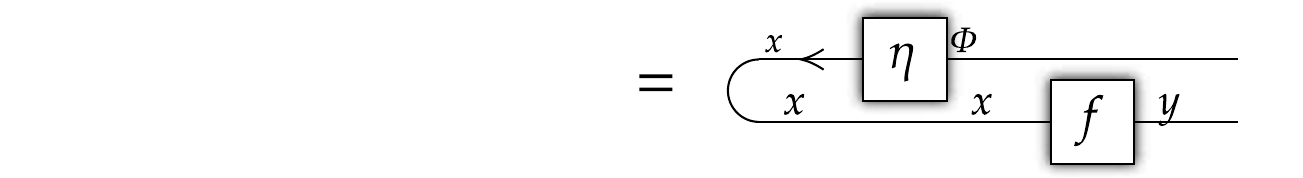}
    \includegraphics[scale=\scalemed]{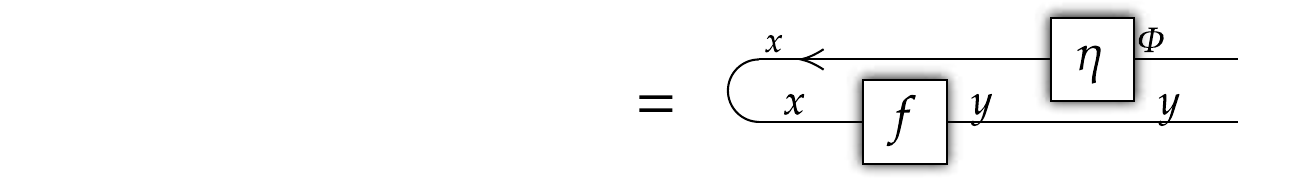}
    \includegraphics[scale=\scalemed]{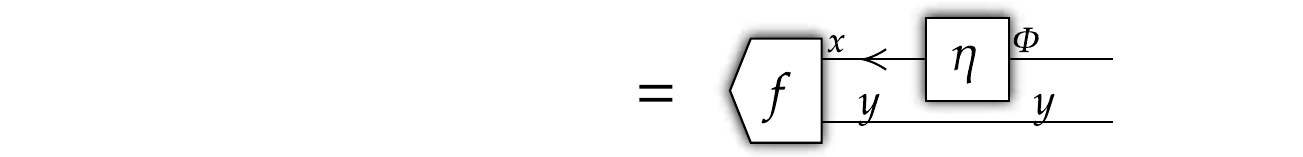}
    \includegraphics[scale=\scalemed]{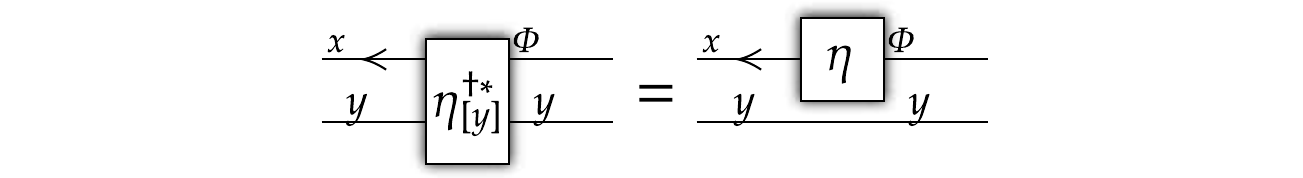}
\end{figure}
\FloatBarrier
This shows that mapping $\displaystyle \eta ^{\dagger }$ back to a natural transformation recovers it. Therefore, the map from elements to natural transformations is a bijection. This implies there is a one-to-one correspondence between natural transformations and elements, completing the proof. $\blacksquare$

\begin{figure}[!htb]
    \centering
    \includegraphics[scale=\scalebig]{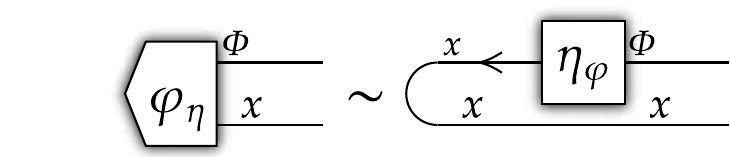}
    \caption{The Yoneda lemma stated with our diagrammatic approach. It shows that elements $\displaystyle \varphi _{\eta } \in \Phi x$ correspond to natural transformations $\displaystyle \eta _{\varphi } :\mathbf{C}( x,\_)\rightarrow \Phi $. This is equivalent to the graphically intuitive rule where elements can ``jump up'' and become natural transformations, and vice-versa.}
    \label{fig:yoneda}
\end{figure}

The Yoneda lemma, fundamentally, says nothing more or less than that the following graphical move shown in Figure \ref{fig:yoneda} is allowed. Elements can become natural transformation components, and vice-versa.
\FloatBarrier
\FloatBarrier

\section{Functor String Diagrams, Products, and Monoidal String Diagrams}\label{sec:products}
In addition to having utility in streamlining algebraic proofs and understanding, functor string diagrams encompass monoidal string diagrams as a specific case. Hence, they are shown to capture existing diagrams and provide a strict improvement over current methods which have found success representing Feynman diagrams, Petri nets, and other structures mentioned in the introduction.

\begin{figure}[h]
  \floatbox[{\capbeside\thisfloatsetup{capbesideposition={left,top}}}]{figure}[\FBwidth]
  {\caption{
		Product categories such as $\displaystyle \mathbf{C} \times \mathbf{D}$ are represented by a double-dashed line. This can represent the combination of expressions from any two categories.
	}
	\label{fig:Product0}}
  {\includegraphics[scale=\scalemed]{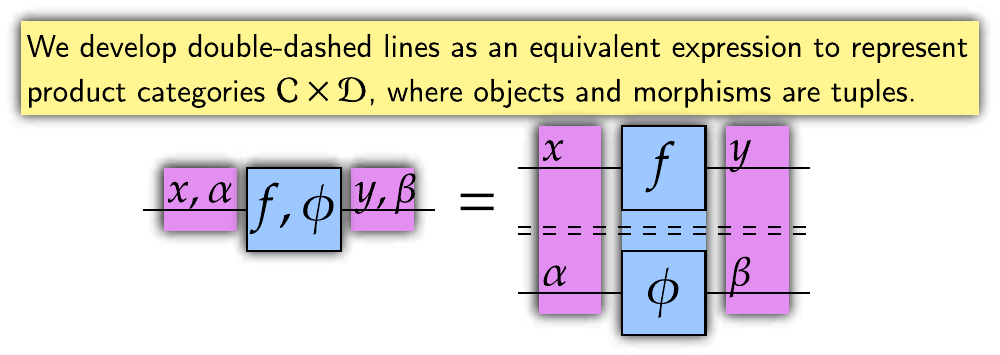}}
\end{figure}

Monoidal string diagrams rely on a category $\displaystyle \mathbf{C}$ having a functor $\displaystyle \otimes $ from $\displaystyle \mathbf{C} \times \mathbf{C}\rightarrow \mathbf{C}$. $\displaystyle \mathbf{C} \times \mathbf{D}$ is the product category wherein objects are tuples and morphisms are tuples from $\displaystyle \mathbf{C}$ and $\displaystyle \mathbf{D}$. We equivalently express objects and morphisms in a product category by a double-dashed line separating the constituent objects or morphisms as in Figure \ref{fig:Product0}.

\FloatBarrier

Therefore, $\displaystyle \mathbf{C} \times \mathbf{C}$ is represented as two expressions from the same category separated by a double dashed line. A functor applied onto a product expression needs to be indicated as distinct from a functor on the upmost expression, and this is achieved with a bold line. Thereby, $\displaystyle \otimes :\mathbf{C} \times \mathbf{C}\rightarrow \mathbf{C}$ can be represented with a bold wire labeled ``$\displaystyle \otimes $'' (see appendix Figure \ref{fig:Product-Options}). 
Furthermore, monoidal categories typically have a privileged product. In this case, it is represented with a single dashed line separating expressions.

\begin{figure}[h]
  \floatbox[{\capbeside\thisfloatsetup{capbesideposition={left,top}}}]{figure}[\FBwidth]
  {\caption{
		By a certain choice of equivalent expression, wherein we choose to view monoidal products as functors and natural transformations, one can see that monoidal string diagrams are a certain specification of functor string diagrams.
	}
	\label{fig:produce-monoidalstring}}
  {\includegraphics[scale=\scalemed]{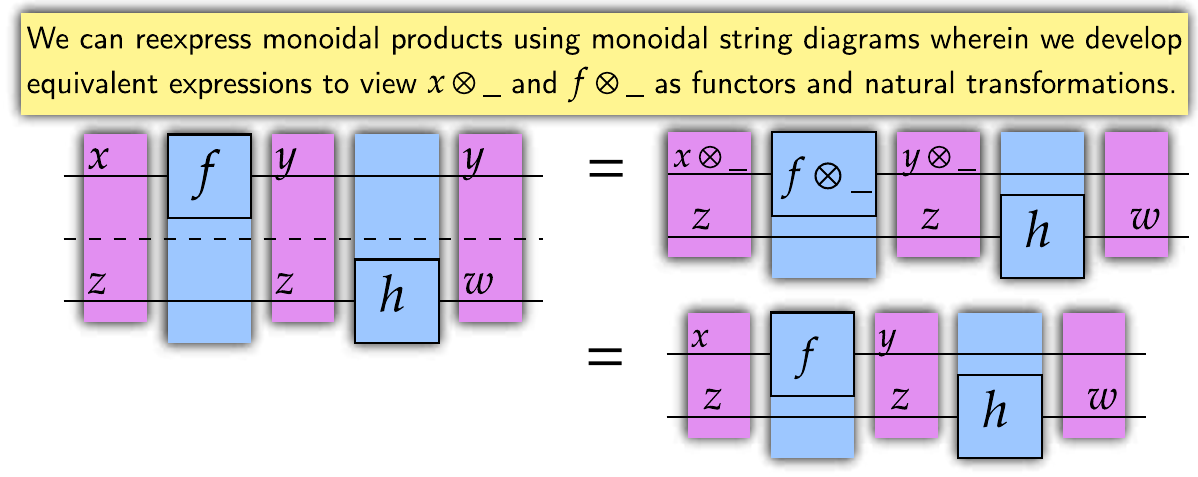}}
\end{figure}

Second, we can represent monoidal categories without a dashed line. This is done by realizing that a monoidal product $\displaystyle \otimes :\mathbf{C} \times \mathbf{C}\rightarrow \mathbf{C}$ grants every object in $\displaystyle x\in \text{Ob}(\mathbf{C})$ a functor $\displaystyle x\otimes \_:\mathbf{C}\rightarrow \mathbf{C}$ and every morphism $\displaystyle f:x\rightarrow y$ a natural transformation $\displaystyle f\otimes \_:x\otimes \_\rightarrow y\otimes \_$. This fact can be observed from the diagram where $\displaystyle f$ can naturally move over $\displaystyle h$, with components $\displaystyle f\otimes \text{Id}_{z}^{z}$ and $\displaystyle f\otimes \text{Id}_{w}^{w}$. This is shown in Figure \ref{fig:produce-monoidalstring}. Monoidal string diagrams, therefore, are functor string diagrams where $\displaystyle x\otimes \_$ is equivalently expressed by a wire labeled $\displaystyle x$, and natural transformations $\displaystyle f\otimes \_$ are expressed with boxes labeled $\displaystyle f$ on functor wires $\displaystyle x\otimes \_$ to $\displaystyle y\otimes \_$.

If we look at how monoidal string diagrams are defined, we see that they obey vertical section decomposition and equivalent expression. A survey by Selinger \cite{selinger-survey-2009} represents a standard introduction. It introduces the features of monoidal string diagrams as re-expressions of symbolic categories, as seen in appendix Figure \ref{fig:Product-Options}.

Monoidal string diagrams focus on powerful algebraic correspondence between diagrammatic and categorical expressions. This \textit{coherence} property states that well-formed diagrams exactly correspond to well-formed expressions. This property is so powerful that it allows all the details of monoidal categories, isomorphisms, associators, and so forth, to follow from the axioms of the graphical language and planar isotopy. Monoidal string diagrams, then, can be used to define what monoidal categories are and hence offer an extremely powerful perspective from which to understand them.

This level of correspondence is not the focus of functor string diagrams. Demanding this isotopy raises a barrier to diagramming and understanding additional structures. To the degree that graphical intuition is useful, it is far simpler to define properties by their local interactions. For instance, natural transformations acting on functor wires follows from a simple local definition. This approach, which does not demand isotopy, allows more features to be introduced with each having some appropriate form of graphical intuition. The graphical intuition embedded into equivalent expressions can coincidentally serve as axioms for the feature, however, this is not the focus of functor string diagrams.

Furthermore, functor string diagrams are sufficiently coherent. All categorical expressions are, ultimately, a chain of morphisms acting on objects, which vertical section decomposition ensures we can perfectly represent. When certain expressions are isomorphic, it is often straightforward to encode that isomorphism as an equivalent expression as we ultimately demand that diagrams correspond to a unique categorical expression up to isomorphism.

\section{Neural Circuit Diagrams}\label{sec:ncd}
Deep learning systems are on the forefront of human technology, with applications in a variety of fields ranging from the scientific investigation of protein folding \cite{senior-improved-2020, jumper-highly-2021} to immense commercial developments. Currently, ad-hoc methods are used to explain these systems, limiting our understanding and obfuscating implementation. Existing categorical machine learning, and theoretical machine learning more generally, often focus on toy models which do not reflect those used in practice given the limit tools to express the details of models \cite{saxe-exact-2014, cockett-reverse-2019, wilson-categories-2022, shiebler-category-2021} (see also appendix Figure \ref{fig:existing-catml}). This stems from a fundamental limitation of monoidal string diagrams, as they are unable to express both the details of separate data and the interaction of axes \cite{chiang-named-2022}. Functor string diagrams, then, have an exciting application regarding formalizing neural circuit diagrams \cite{abbott-neural-2023, abbott-robust-2023}, a categorical method for comprehensively communicating deep learning models.

\begin{figure}[h]
  \floatbox[{\capbeside\thisfloatsetup{capbesideposition={left,top}}}]{figure}[\FBwidth]
  {\caption{
		For deep learning models, objects typically consist of hom-functors to $\mathbb{R}$. We develop an equivalent expression which lets this be readily expressed.
	}
	\label{fig:DL-form}}
  {\includegraphics[scale=\scalemed]{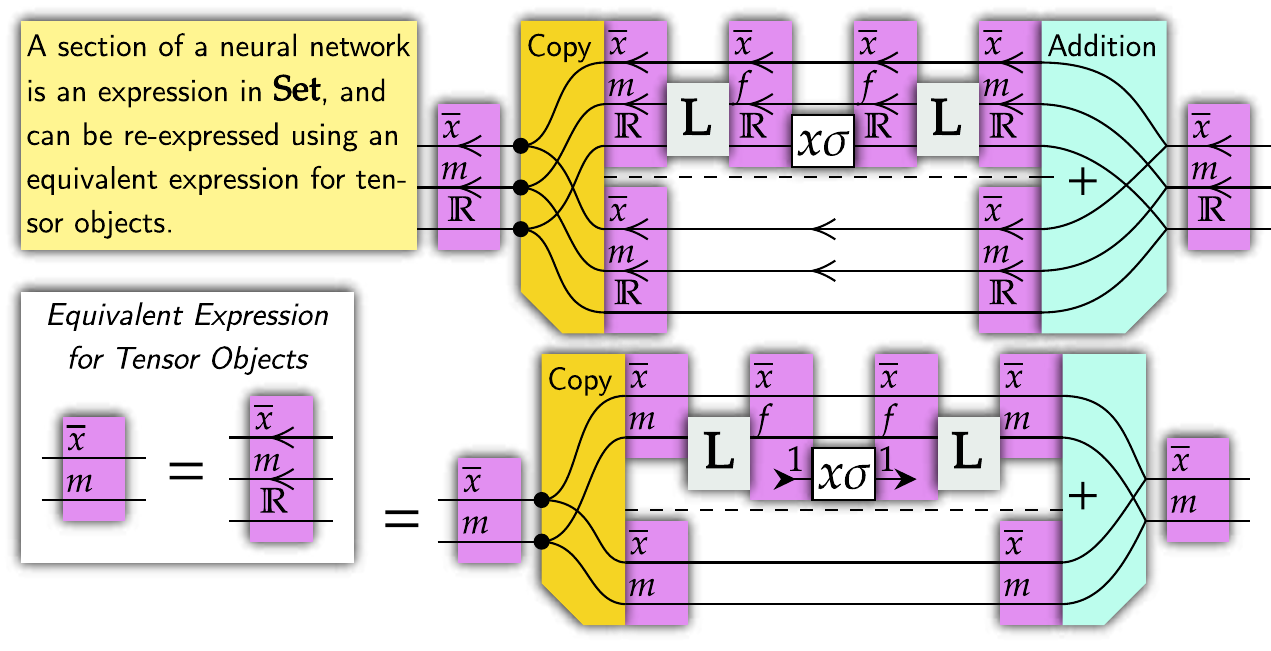}}
\end{figure}

\paragraph{The Form of Deep Learning Models.} Deep learning models are functions between data types, meaning they are expressions in $\textbf{Set}$. Typically, data types are tensors of an underlying data type $\mathbb{R}$. Therefore, objects are almost always of hom-functors applied to $\mathbb{R}$. We develop an equivalent expression, whereby wires withour arrows represent those objects applied as hom-functors to $\mathbb{R}$. This allows broadcasted functions between tensors to be represented as in Figure \ref{fig:DL-form}. Furthermore, to represent functions $\mathbb{R}\rightarrow \mathbb{R}$, we exploit that the $1$-hom-functor is an isomorphism, letting it be freely introduced and removed (unsqueezing and squeezing), which we introduce and remove with a wire with an arrowhead.

\paragraph{Understanding Expressions with Indexes. }As we are operating in $\displaystyle \mathbf{Set}$, every object $\displaystyle a$ has a set of coprojections given by its elements $\displaystyle \bra{i_{a}} :1\rightarrow a$. Morphisms $\displaystyle a\rightarrow b$ can be uniquely identified and freely constructed from pre-compositions with elements. By the Yoneda lemma, natural transformations $\displaystyle \mathbf{Set}( a,\_)\rightarrow \mathbf{Set}( 1,\_)$ correspond to elements of $\displaystyle \mathbf{Set}( 1,a)$. Therefore, elements $\displaystyle \bra{i_{a}} :1\rightarrow a$ become natural transformations $\displaystyle | i^{\overleftarrow{a}}\rangle^{x}_x :\mathbf{Set}( a,x)\rightarrow \mathbf{Set}( 1,x)$ which we call \textit{indexes.} They form a set of projections (see appendix Figure \ref{fig:indexes1} and \ref{fig:indexes2}), letting morphisms to $\textbf{Set}(a,x)$ be uniquely identified and freely constructed from post-composition with indexes.

\begin{figure}[htb]
  \floatbox[{\capbeside\thisfloatsetup{capbesideposition={left,top}}}]{figure}[\FBwidth]
  {\caption{
		Continuous wires implement naturality, letting natural transformations, including indexes, be transmitted, defining these expressions. 
	}
	\label{fig:DL-broadcast}}
  {\includegraphics[scale=\scalemed]{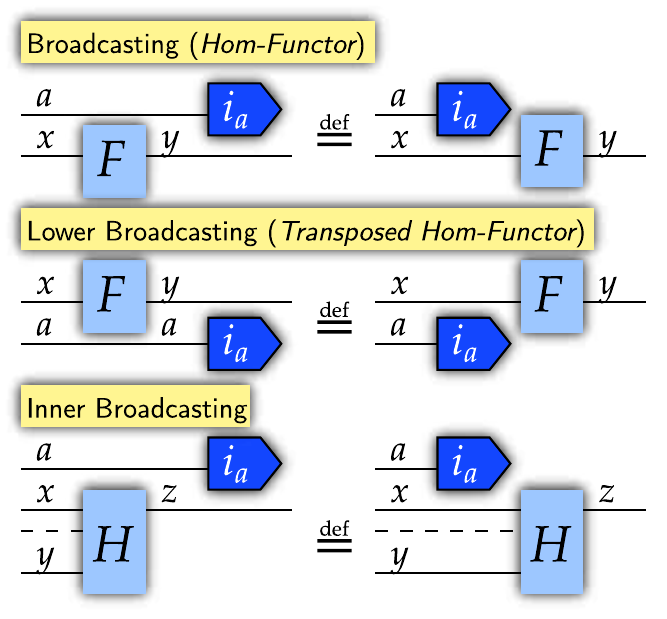}}
\end{figure}

\paragraph{Broadcasting.} Using indexes, we can define and construct various useful forms of functions. This includes broadcasting, corresponding to hom-functors, and inner broadcasting, corresponding to broadcasting within tuple segments. These expressions use continuous wires to indicate the naturality of indexes, as in Figure \ref{fig:DL-broadcast}.

\begin{figure}[!h]
  \floatbox[{\capbeside\thisfloatsetup{capbesideposition={left,top}}}]{figure}[\FBwidth]
  {\caption{
		Natural transformations can be shown with hexagons, and they have naturality regarding hom-functor broadcasting.
	}
	\label{fig:DL-nat}}
  {\includegraphics[scale=\scalemed]{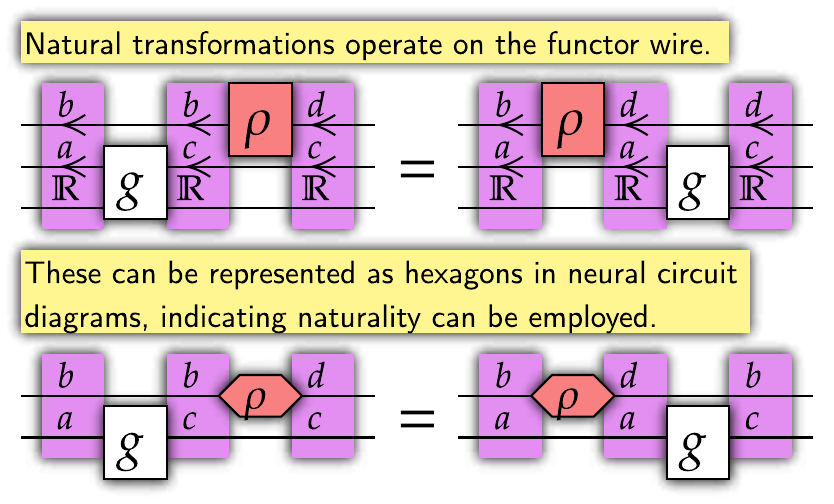}}
\end{figure}

\begin{figure}[h]
  \floatbox[{\capbeside\thisfloatsetup{capbesideposition={left,top}}}]{figure}[\FBwidth]
  {\caption{
		By testing this inner broadcasting expression with indexes $k_b$, we see that the naturality of indexes as indicated in Figure \ref{fig:DL-broadcast} implies a naturality of natural transformations along continuous wires.
	}
	\label{fig:DL-cont}}
  {\includegraphics[scale=\scalemed]{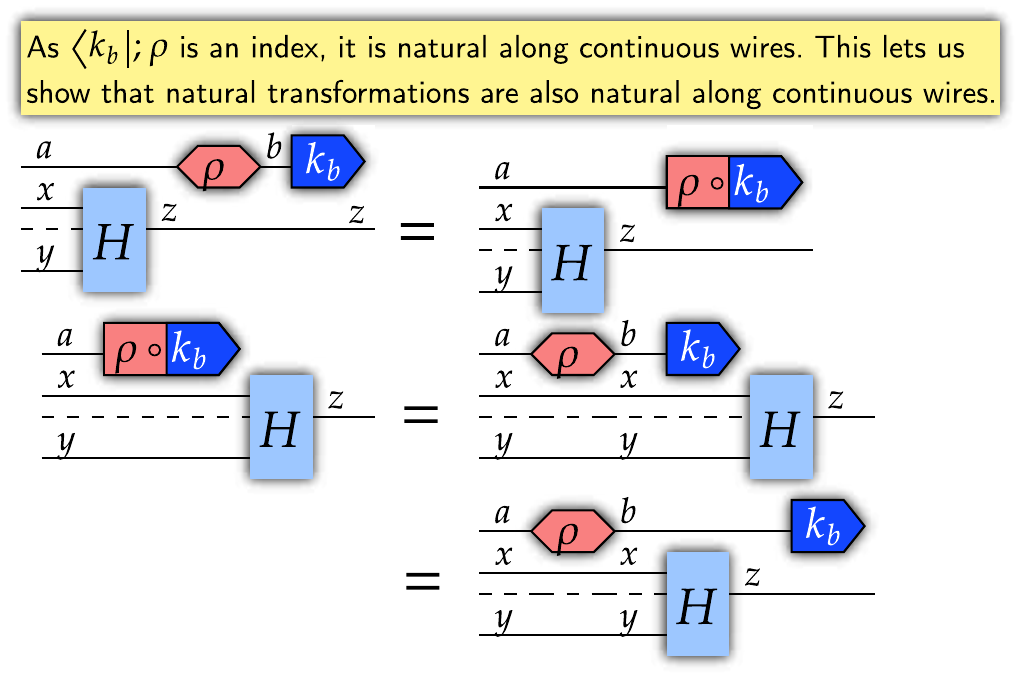}}
\end{figure}

\paragraph{Natural Transformations.} By the Yoneda lemma, functions on sets $x \rightarrow y$ correspond to functions manipulating the indexes of tensors $\mathbb{R}^y\rightarrow \mathbb{R}^x$. We represent natural transformations in neural circuit diagrams using hexagons (see Figure \ref{fig:DL-nat}), and note that they travel along continuous wires the same way as indexes (see Figure \ref{fig:DL-cont}). Continuous wires defined by naturality with indexes, then, can be seen to grant naturality to natural transformations more generally.

\paragraph{Einops.} Einstein operations refers to a common class of rearranging and contracting (summation) over axes. This covers a common class of operations, the details of which are critical yet challenging to express with traditional methods. Einstein operations can be expressed by associating desired input with output indexes symbolically, such as with Einstein tensor notation or the Einops package \cite{rogozhnikov-einops-2021}. With neural circuit diagrams, they correspond to continuous wiring associating input and output axes. These operations are linear and often require an outer product to first be taken (element-wise multiplication), which we diagram by terminating a single dotted line.

\begin{figure}[!htb]
    \centering
    \includegraphics[scale=\scalemed]{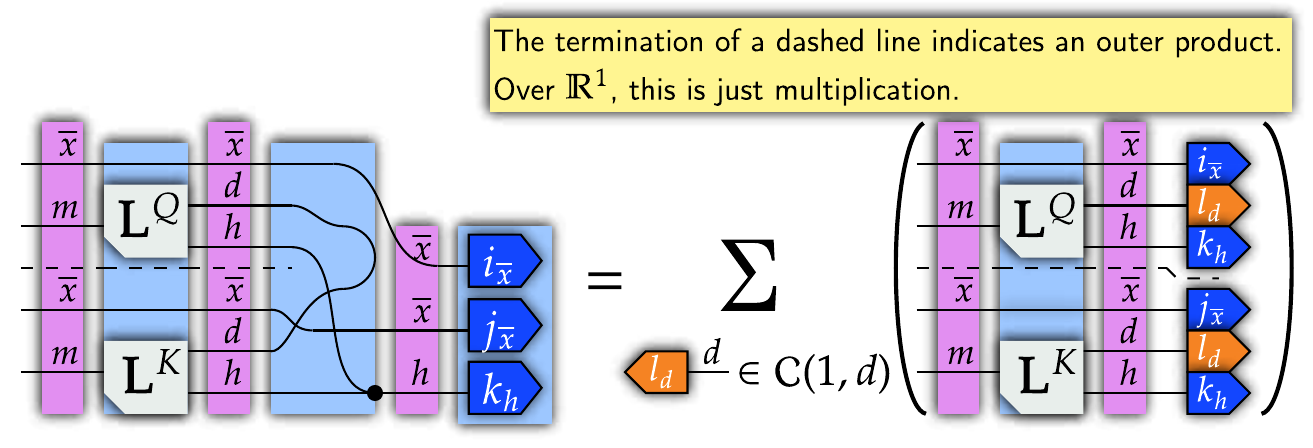}
    \caption{Natural operations, including indexes, travel along continuous wires. This means the above expression, for some choice of indexes $i_{\bar{x}}$, $j_{\bar{x}}$, $k_h$ placed after the Einstein operation, can be seen to be a contraction over the indexes placed before the operation.}
    \label{fig:DL-einops}
\end{figure}

\begin{figure}
    \centering
    \includegraphics[scale=\scalemed]{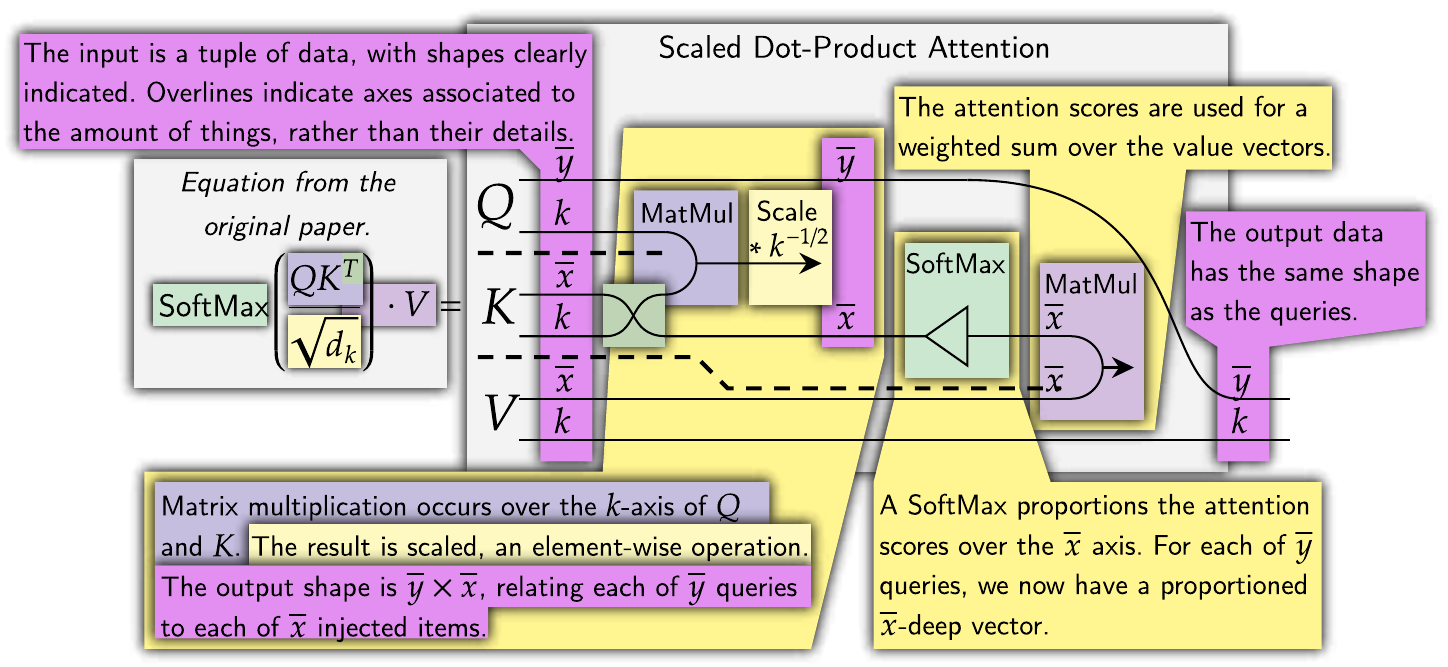}
    \caption{Scaled dot-product attention can be far more clearly expressed with neural circuit diagrams. The interaction between axes and the changing shape of data is naturally shown by diagrams.}
    \label{fig:SDPA}
\end{figure}

\begin{figure}[!htb]
    \centering
    \includegraphics[scale=\scalemed]{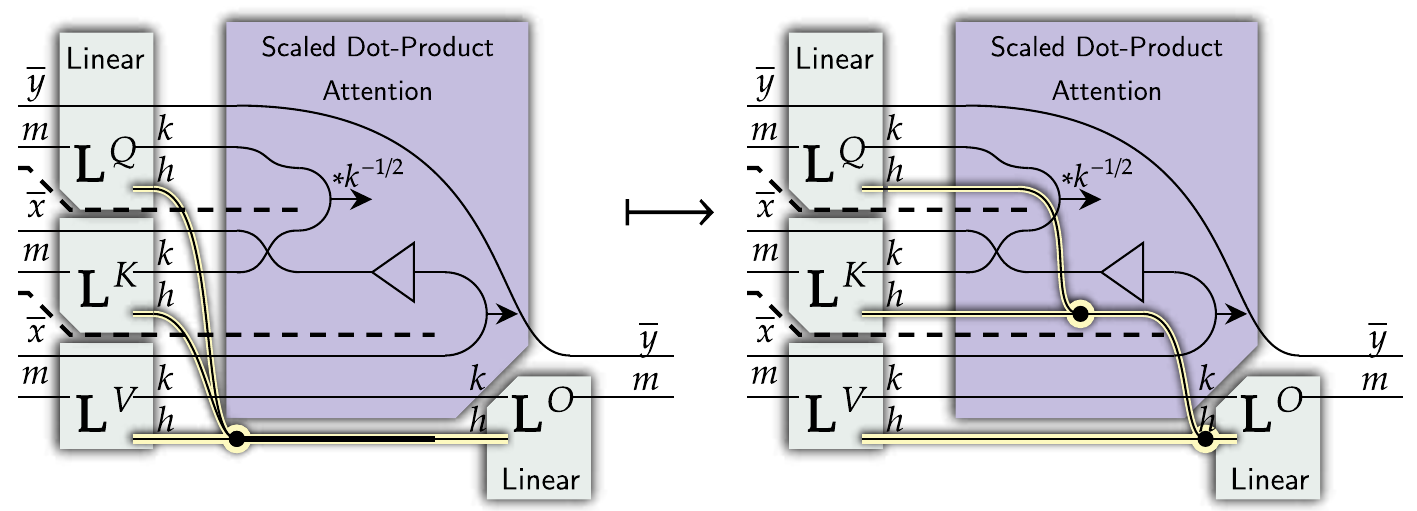}
    \caption{Conceptually, multi-head attention involved broadcasting SDPA over $\displaystyle h$ parallel operations and aggregating the result. The queries, keys, and valued are produced and recombined by adding an additional axis to the learned linear layers.}
    \label{fig:MHA}
\end{figure}

\paragraph{Scaled dot-product attention. }Using these tools, we can express a significant yet complicated deep learning algorithm, multi-head attention, clearly. This algorithm requires interaction between specific axes and independent data, which is hard to express with traditional methods (see appendix Figure \ref{fig:typical} for a typical presentation). With neural circuit diagrams, the specific interactions for scaled dot-product attention (SDPA) can be clearly seen as in Figure \ref{fig:SDPA}. Furthermore, multi-head attention which performs SDPA $h$-times and aggregates the result, can be shown be shown with neural circuit diagrams in Figure \ref{fig:MHA}.

\FloatBarrier
\section{Conclusion}
In this work, we stated the need for a general means of diagrammatically expressing categorical ideas. Functor string diagrams are a general and principled approach to developing categorical diagrams that have robust mathematics for theoretical investigations, encompass existing monoidal string diagram, and have clear utility for application. Compared to the existing paradigm of monoidal string diagrams and string diagrams, functor string diagrams offer a flexible yet principled approach that ensures complex ideas can be clearly expressed. Their ability to express both functors and products allows for theoretical ideas to be clearly understood, and allows for applied systems to be straightforwardly modeled.

\bibliographystyle{eptcs}
\bibliography{imported}

\appendix
\FloatBarrier
\section*{Appendix}
\section{Current Categorical Diagrams} \label{sec:existing}
\begin{figure}[htb]
    \centering
    \includegraphics[width=0.75\linewidth]{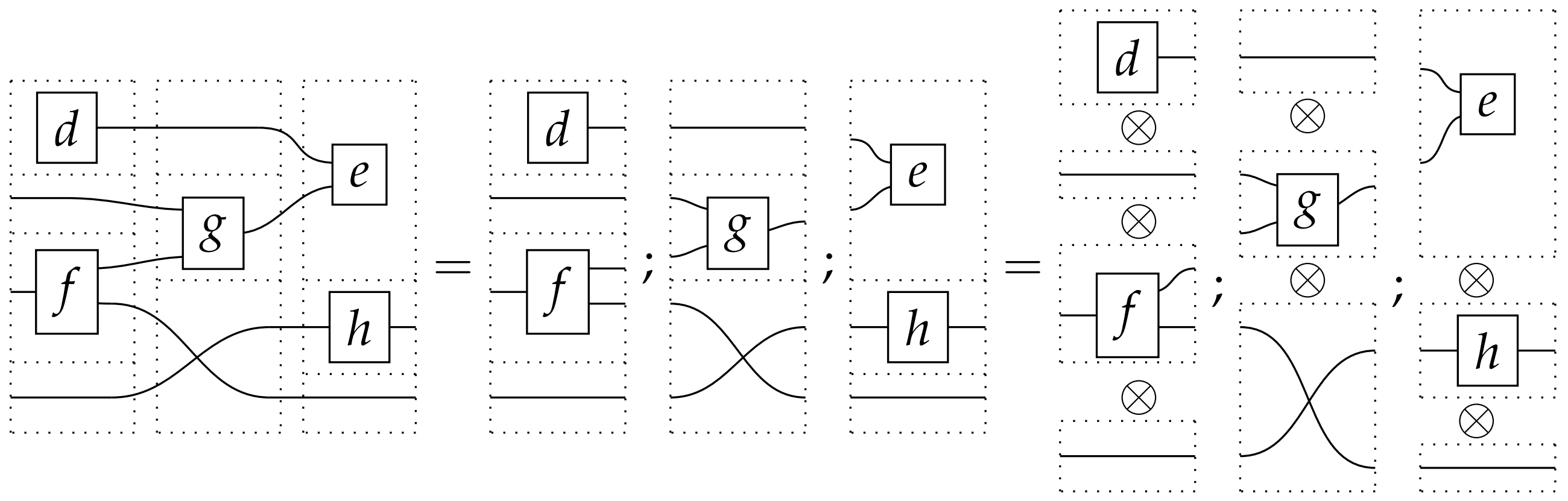}
    \caption{A traditional monoidal string diagram from page 18 of \cite{piedeleu-introduction-2023}. This expression clearly shows the objects and morphisms at each stage, and has graphical isotopy whereby topological transforms of diagrams lead to valid, isomorphic expressions. However, it does not show the details of axes should the objects be Cartesian tuples. Furthermore, it cannot show functors which are not monoidal products of objects nor natural transformations which are not morphisms between the monoidal product of objects. The useful properties of monoidal string diagrams are encompassed by functor string diagrams, as shown in Section \ref{sec:products}.}
    
\end{figure}
\begin{figure}
    \centering
    \includegraphics[width=0.5\linewidth]{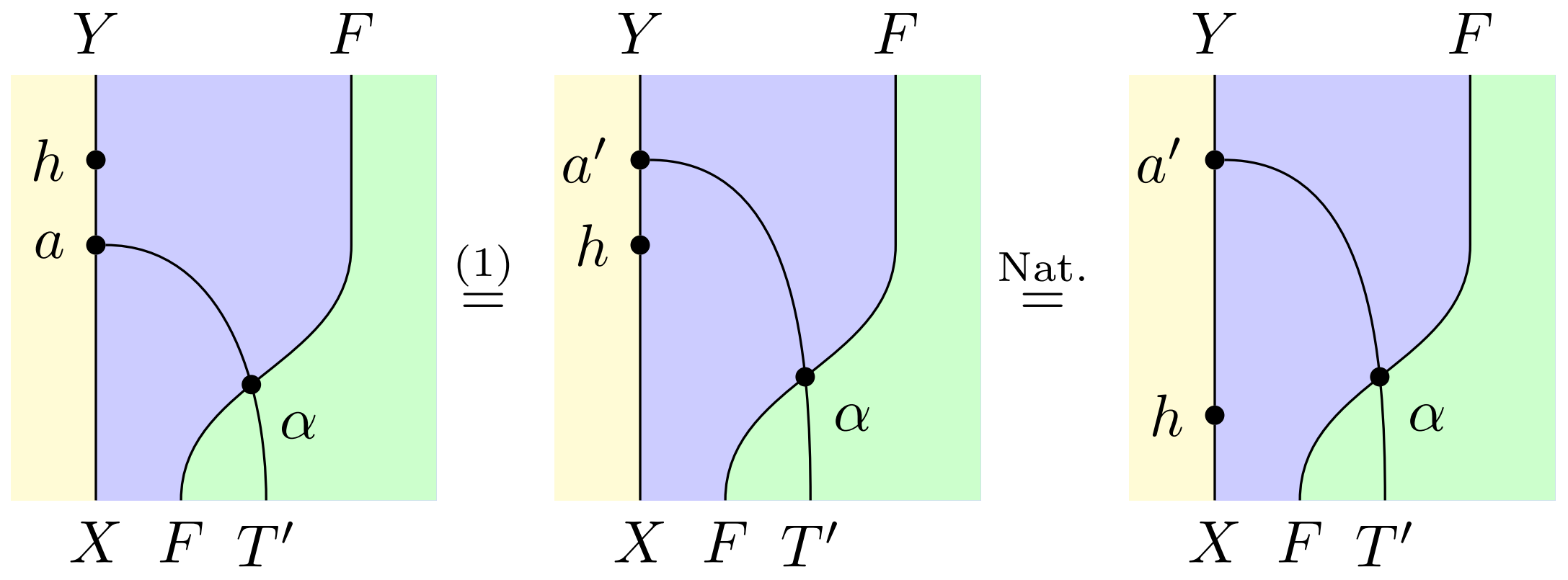}
    \caption{A string diagram from Marsden \cite{marsden-category-2014}. These expressions use coloured regions to show the current category, and display functors as left-to-right boundaries between coloured regions. Natural transformations are dots on functors or their intersection. These can clearly display the action of functors, natural transformations, and more complicated structures such as adjunctions, monads, and (co)limits, and display graphical isotopy, implementing the axioms of these structures graphically. However, they fail to be clearly decomposable for applied cases, and have a cumbersome representation of products. The behaviour of functors and natural transformations using functor string diagrams are shown in Section \ref{sec:foundations}.}
    
\end{figure}
\begin{figure}
    \centering
    \includegraphics[width=0.75\linewidth]{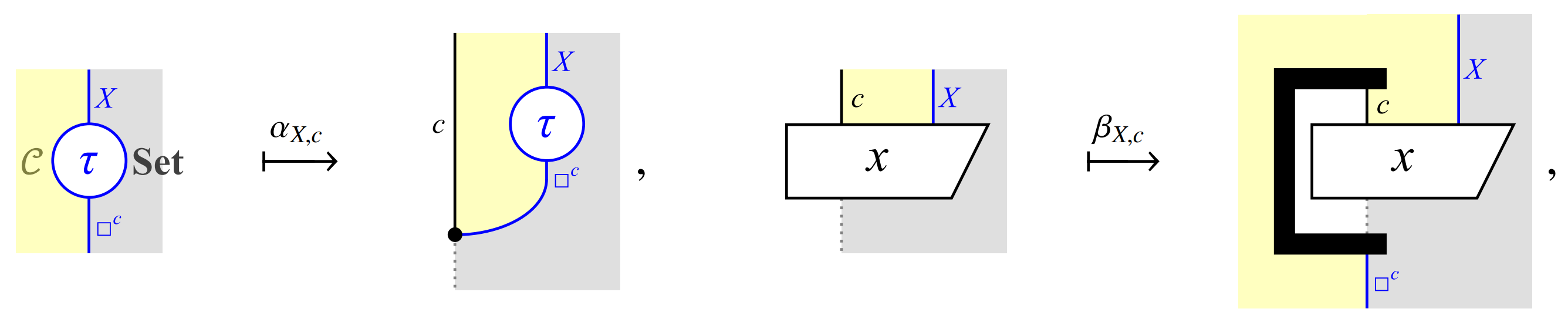}
    \caption{Nakahira's diagrams \cite{nakahira-diagrammatic-2023} extends Marsden's diagrams to consider substiutions and various algebraic proofs, including the Yoneda lemma shown above. These diagrams have powerful tools for the consideration of products, bifunctors, and other constructs. However, streamlined simplicity is preferred for practical application and intuitive understanding, which is shown in the functor string diagram proof for the Yoneda lemma in Section \ref{sec:Yoneda}.}
    
\end{figure}
\begin{figure}
    \centering
    \includegraphics[width=0.5\linewidth]{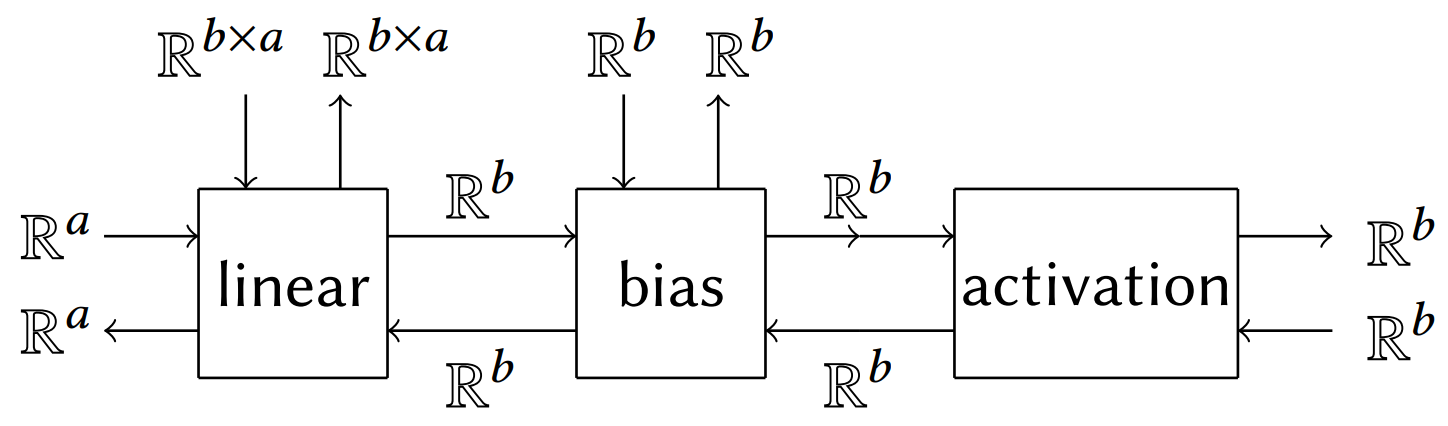}
    \caption{Categorical machine learning such as from \cite{cruttwell-categorical-2021} uses monoidal string diagrams with Cartesian products, along with the $\textbf{Para}$ construct which allows parameters to be introduced from above. These do not express the details of axes, and hence focus on toy models rather than those used in practice. Since 2015 \cite{he-deep-2015, he-identity-2016}, for example, deep learning has employed residual connections which are significantly different from linear-bias-activation toy models. Functor string diagrams underly neural circuit diagrams (see Section \ref{sec:ncd}), which can rigorously and comprehensively express deep learning models.}
    \label{fig:existing-catml}
\end{figure}
\FloatBarrier
\section{Equivalent Expression of Products}
\begin{figure}[!htb]
    \centering
    \includegraphics[scale=\scalemed]{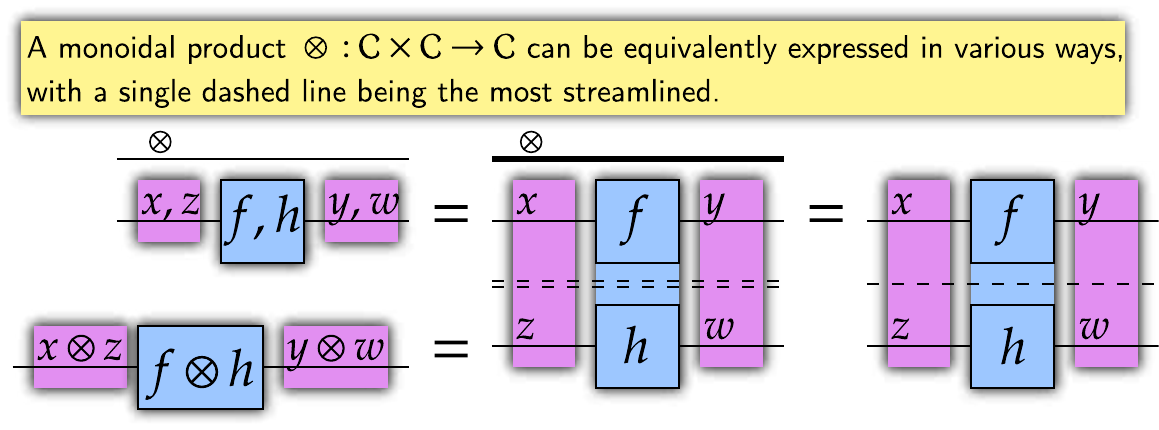}
    \caption{These equivalent expressions for the diagramming of product categories and monoidal products offer immense utility. We see how different expressions are related and can freely adopt whichever approach is the most useful.}
    \label{fig:Product-Options}
\end{figure}

With functor string diagrams, we can express bifunctors $\textbf{C}\times \textbf{C} \rightarrow \textbf{C}$ as bolded functor wires over an expression in $\textbf{C} \times \textbf{C}$, as in Figure \ref{fig:Product-Options}. As we can see, functor string diagrams provide a flexible framework to understand the connections between various means of viewing expressions.

\FloatBarrier
\section{Typical Introduction to Monoidal String Diagrams}
\begin{figure}[h]
  \floatbox[{\capbeside\thisfloatsetup{capbesideposition={left,top}}}]{figure}[\FBwidth]
  {\caption{
		Table 1 from \cite{selinger-survey-2009}, ``The graphical language of categories''.
	}
	\label{fig:table1}}
  {\includegraphics[scale=0.125]{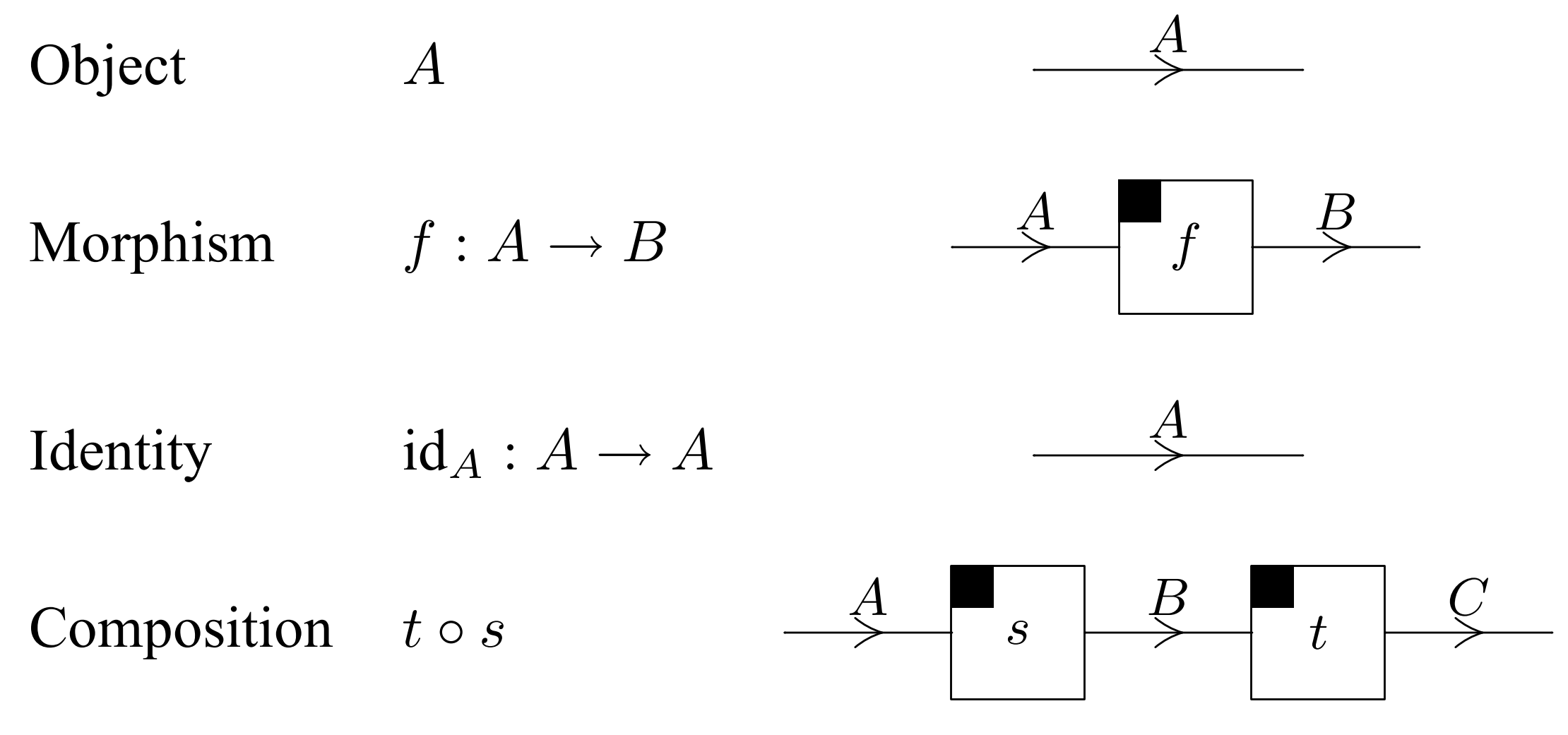}}
\end{figure}

\begin{figure}[h]
  \floatbox[{\capbeside\thisfloatsetup{capbesideposition={left,top}}}]{figure}[\FBwidth]
  {\caption{
		Table 2 from \cite{selinger-survey-2009}, ``The graphical language of monoidal categories''.
	}
	\label{fig:table2}}
  {\includegraphics[scale=0.125]{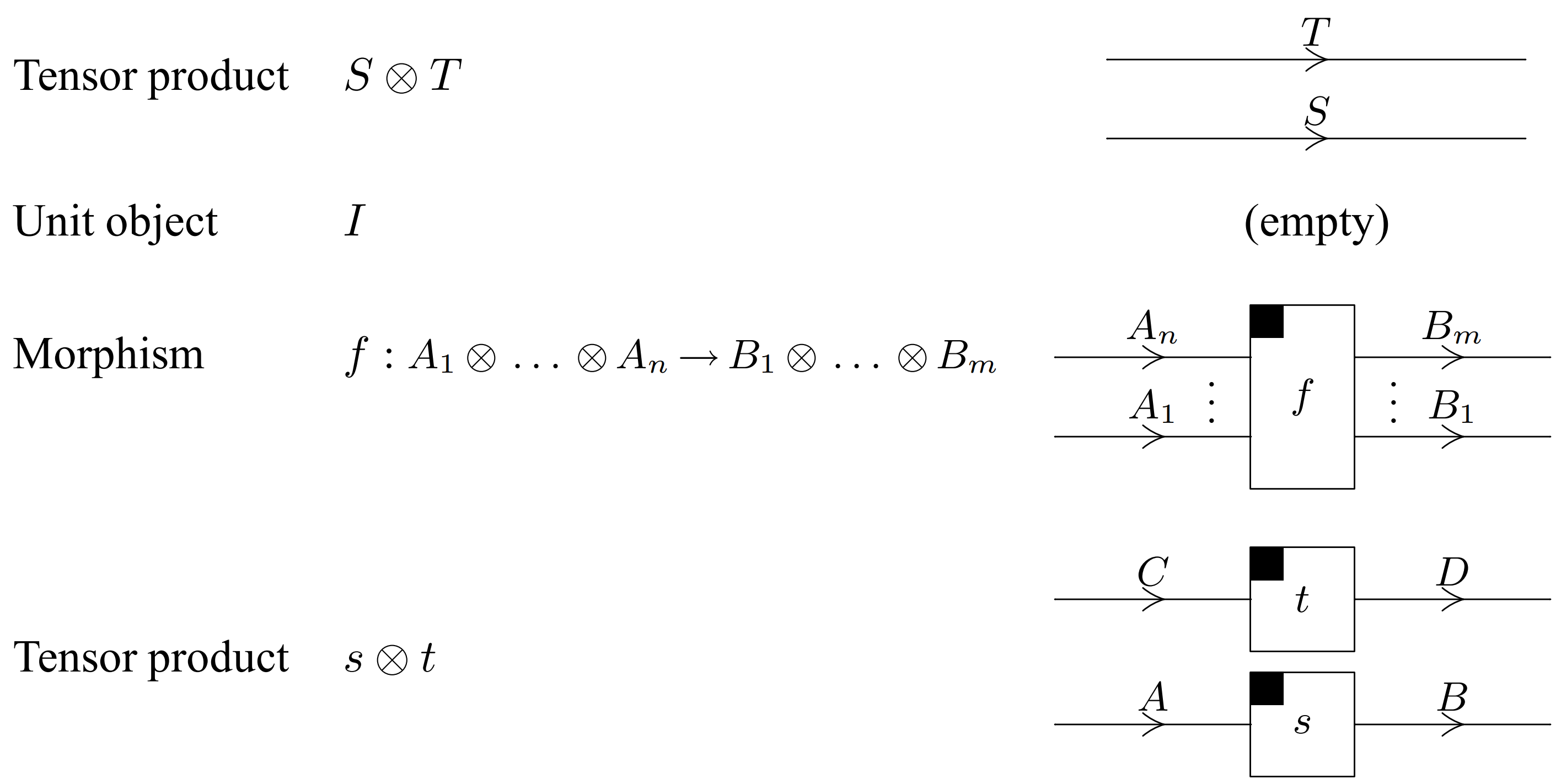}}
\end{figure}

Typical introductions to monoidal string diagrams introduce the fundamental structures by elaborating how objects and morphisms are represented. Below, tables 1 and 2 from \cite{selinger-survey-2009} are shown in Figure \ref{fig:table1} and \ref{fig:table2}. This lets us see that monoidal string diagrams are, fundamentally, diagrams with vertical sections corresponding to objects and morphisms equipped with equivalent expressions.

\FloatBarrier
\section{Elements, Indexes, and Projections}
\begin{figure}[htb]
  \floatbox[{\capbeside\thisfloatsetup{capbesideposition={left,top}}}]{figure}[\FBwidth]
  {\caption{
		Elements of a set $1 \rightarrow a$ form a set of coprojections which can uniquely identify and construct functions $a\rightarrow \_$. 
	}
	\label{fig:indexes1}}
  {\includegraphics[scale=0.5]{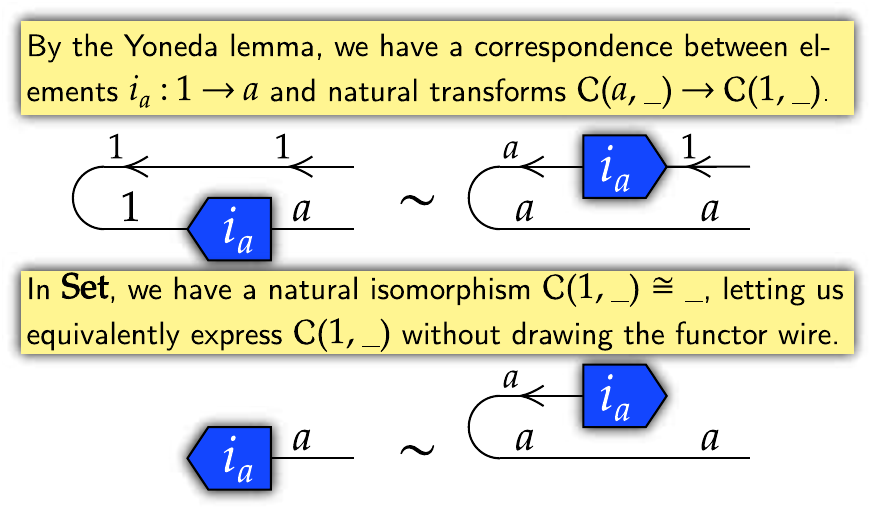}}
\end{figure}

\begin{figure}[htb]
  \floatbox[{\capbeside\thisfloatsetup{capbesideposition={left,top}}}]{figure}[\FBwidth]
  {\caption{
		The associated natural transformations $\textbf{Set}(a,\_)\rightarrow \_$, for a set of projections, letting functions $\_ \rightarrow \textbf{Set}(a,\_)$ be uniquely identified.
	}
	\label{fig:indexes2}}
  {\includegraphics[scale=0.5]{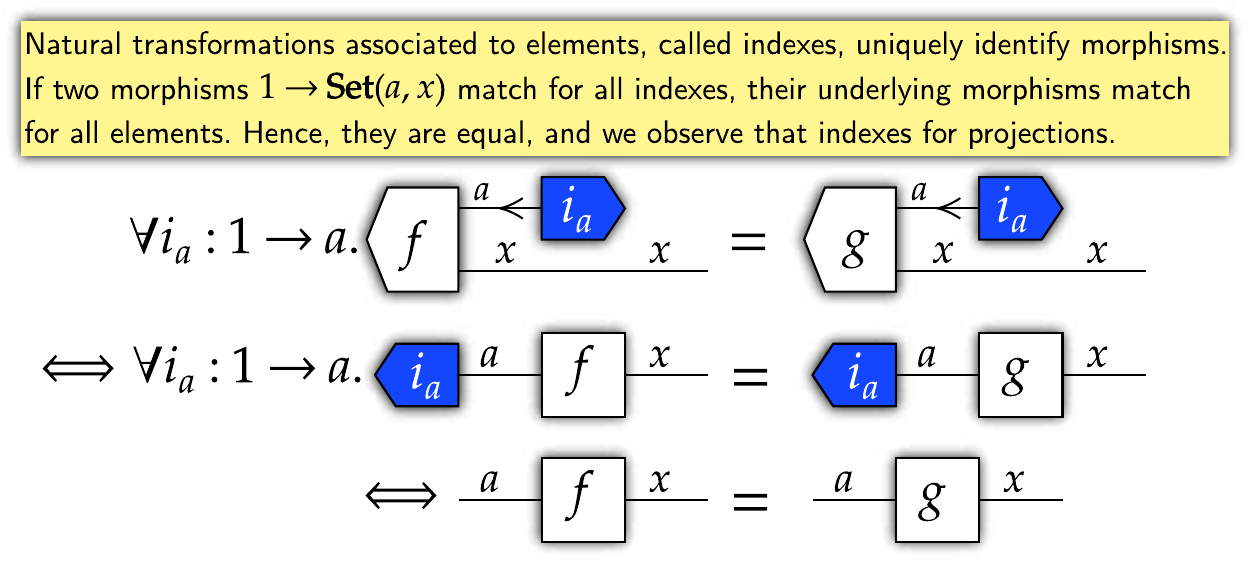}}
\end{figure}

\FloatBarrier
\newpage
\section{Typical Presentation of Attention}
\begin{figure}[h]
\begin{minipage}[t]{0.4\textwidth}
  \centering
  Scaled Dot-Product Attention \\
  \vspace{0.5cm}
  \includegraphics[scale=0.6]{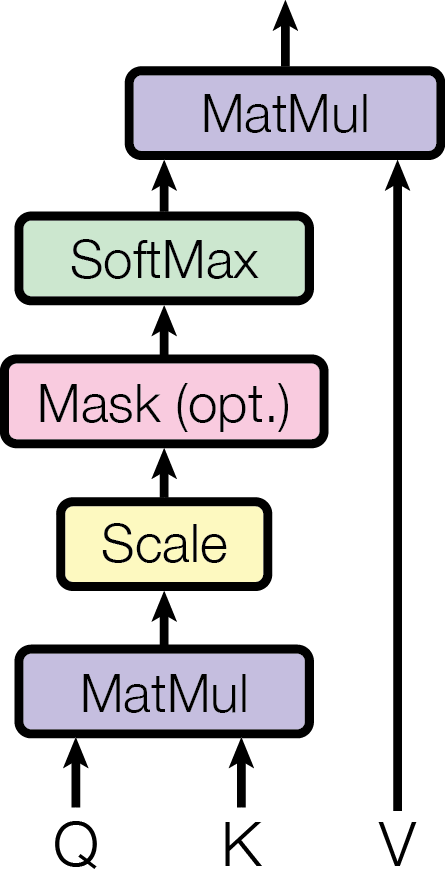}
\end{minipage}
\begin{minipage}[t]{0.4\textwidth}
  \centering 
  Multi-Head Attention \\
  \vspace{0.1cm}
  \includegraphics[scale=0.6]{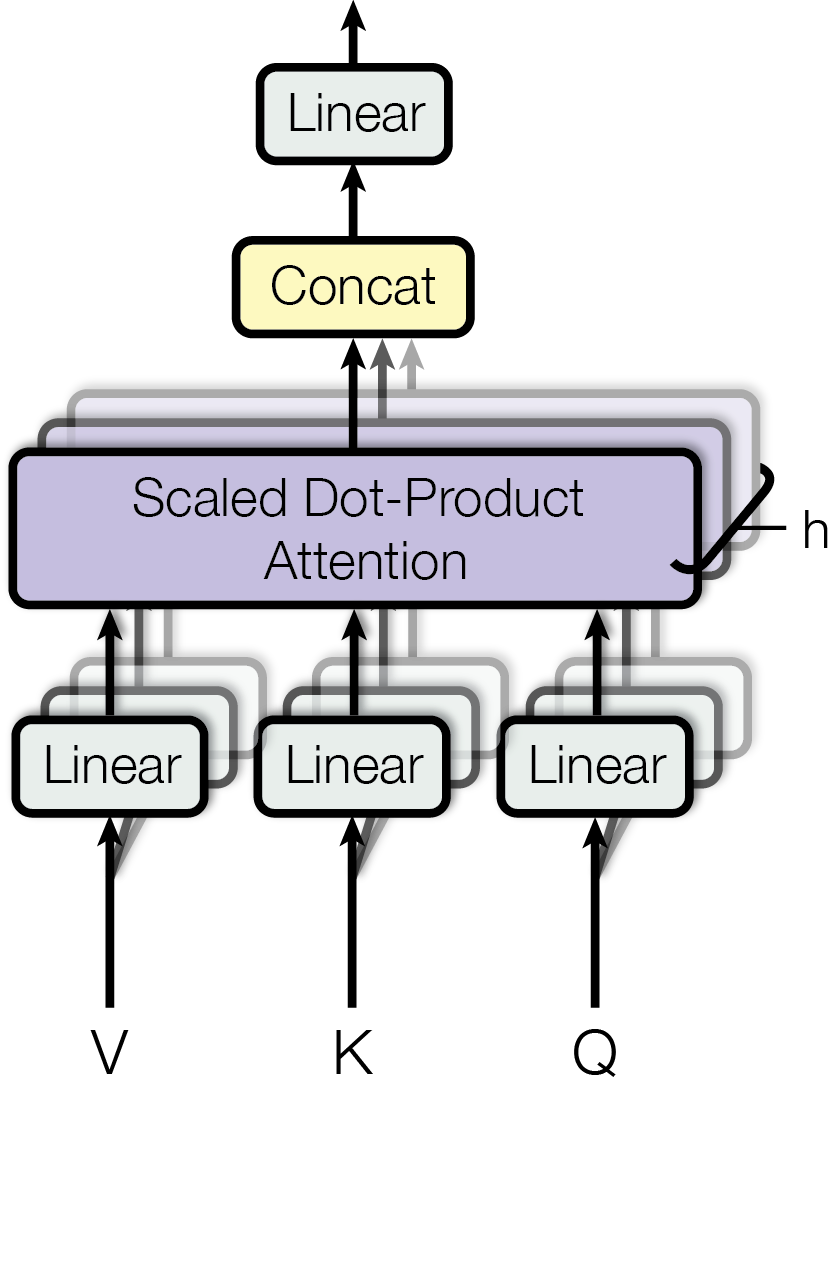}  
\end{minipage}

  \caption{The prevailing ad-hoc diagrams used to represent scaled dot-product attention and multi-head attention obfuscate the interaction of axes and how the shape of data changes across operations, making understanding and implementation difficult.}
  \label{fig:typical}
\end{figure}

A typical presentation of scaled dot-product attention can be drawn from \cite{NIPS2017-3f5ee243} are shown in Figure \ref{fig:typical}. This work introduced transformer models. However, its diagrams are difficult to understand, especially with regards to the interaction of axes and the shape of data throughout the model. As can be seen, Figure \ref{fig:SDPA} and \ref{fig:MHA} are a significant improvement.

\end{document}